\documentclass[11pt
]{amsart}

\usepackage{hyperref}
\hypersetup{bookmarksdepth=3}

\usepackage{geometry}
\usepackage{amsmath,amssymb}
\usepackage{
mathrsfs}
\usepackage{esint}

\usepackage{tensor}

\usepackage{amssymb,amsmath,amsthm}

\newtheorem{theorem}{Theorem}

\newtheorem{definition}[theorem]{Definition}

\newtheorem{remark}[theorem]{Remark}

\usepackage{graphicx}

\begin{document}

\title{Some Remarks on Extrapolation with ``flat" weights}
\author{Nicholas Boros, Nikolaos Pattakos and Alexander Volberg} \address{Department of Mathematics, Michigan State University, East
Lansing, MI 48824, USA}

\subjclass{30E20, 47B37, 47B40, 30D55.} 
\keywords{Key words: Calder\'on--Zygmund operators, $A_2$ weights, extrapolation.}
\date{}

\begin{abstract}
{We prove an extrapolation result for general operators under some weak assumptions on the boundedness of the operator. In particular, we show that if the operator is weakly bounded on some $L^{p_{0}}(w)$, for all ``flat" weights, $w\in A_{p_{0}}$, $1<p_{0}<\infty$, then for $p$ in some small neighborhood around $p_{0}$, and all ``flat" $A_{p}$ weights, w, the operator is weakly bounded on $L^{p}(w)$, and as a result we get strong type estimates for the operator. This comes in comparison with the general extrapolation result of Rubio de Francia which can be found in \cite{GCRF}.   }
\end{abstract}

\maketitle

\begin{section}{introduction, notation and the main result}
\label{sec1}

In recent years weighted estimates for singular integrals have been a very active area of interest. By a weight we are referring to a positive $L^{1}_{loc}(\mathbb R^{n})$ function, w. Under the assumption that the quantity

$$[w]_{A_{p}}:=\sup_{Q}\Big(\frac1{|Q|}\int_{Q}w(x)dx\Big)\Big(\frac1{|Q|}\int_{Q}w(x)^{-\frac1{p-1}}dx\Big)^{p-1},$$
is finite, where $Q$ is a cube in $\mathbb R^{n}$, we say that $w$ is an $A_{p}$ weight, $1<p<\infty$, and the number $[w]_{A_{p}}$ is called the $A_{p}$ characteristic of $w$ (for more information about weights see \cite{GCRF}). Let us give the following definition. 

\begin{definition}
An $A_{p}$ weight $w$ is called a ``flat" weight, if it's characteristic is close to $1$.
\end{definition}

We say that an operator $T$ is weakly bounded on $L^{p}(w)$ (or of weak type $(p,p)$ with respect to $w$), if the following estimate holds

$$\sup_{\lambda>0}\lambda (w\{x\in\mathbb R^n:|Tf(x)|>\lambda\})^{\frac1{p}}\leq c\|f\|_{L^{p}(w)},$$
for all functions $f\in L^{p}(w)$, where $c>0$ does not depend on $f$. The smallest such $c$ is denoted by $\|T\|_{L^{p}(w)\rightarrow L^{p,\infty}(w)}$. Similarly, we say that the operator $T$ is strongly bounded on $L^{p}(w)$ (or just bounded), if there exists $c>0$ independent of $f\in L^{p}(w)$ such that

$$\|Tf\|_{L^{p}(w)}\leq c\|f\|_{L^{p}(w)},$$
$\forall f\in L^{p}(w)$. The smallest such $c$ is denoted by $\|T\|_{p,w}$ or by $\|T\|_{L^{p}(w)\rightarrow L^{p}(w)}$. It is obvious that strong type estimates imply weak type estimates. The converse is not true in general.

It has been known that classical singular integral operators of Calder\'on-Zygmund type are bounded on $L^{p}(w)$ for weights $w\in A_{p}$. Recently though, the focus has been on a problem, now known as the $A_{2}$ conjecture for Calder\'on-Zygmund operators, on how to find the sharp dependence of the operator norm and the $A_{2}$ characteristic of the weight. It states that for any singular integral operator $T$ of Calder\'on-Zygmund type we have the estimate

$$\|T\|_{L^{2}(w)\rightarrow L^{2}(w)}\leq c[w]_{A_{2}},$$
for all $A_{2}$ weights $w$, where $c$ is a positive constant independent of the weight. This result turns out to be correct and was first proven in \cite{H}. This linear estimate with respect to the $A_{2}$ characteristic of the weight, is sharp for many of the classical operators, such as the Hilbert, Martingale, Riesz and Beurling transforms.  In order to prove a similar estimate for the $L^{p}(w)$ spaces, where $p\neq2$, we have to use the following extrapolation result. This formulation, that we will follow, appears in \cite{DMP}.

\begin{theorem}
\label{result1}
If for some $1<p_{0}<\infty$, there exists $\alpha(p_{0})>0$ such that

$$\|T\|_{L^{p_{0}}(w)\rightarrow L^{p_{0}}(w)}\leq c[w]_{A_{p_{0}}}^{\alpha(p_{0})},$$
for all $A_{p_{0}}$ weights, w, where $c>0$ does not depend on $w$, then for all $1<p<\infty$,

$$\|T\|_{L^{p}(w)\rightarrow L^{p}(w)}\leq c[w]_{A_{p}}^{\alpha(p_{0})\max\{1,\frac{p_{0}-1}{p-1}\}},$$
for all $A_{p}$ weights, $w$, where $c>0$ does not depend on the weight.
\end{theorem}

This theorem can be found in \cite{DMP}, and is sharp with respect to the exponent of the $A_{p}$ characteristic for many classical Calder\'on-Zygmund operators. It can be applied to all operators that satisfy a certain sharp weighted estimate for all weights in some of the $A_{p}$ classes. 

Let us also mention that weighted estimates of singular integrals of Calder\'on-Zygmund type are a natural problem to consider in the study of PDE's.  In fact, Calder\'on-Zygmund operators naturally arise as fraction derivatives of solutions of PDE's.  If we recall, for example, in \cite{PV} the authors proved that the Ahlfors-Beurling operator in the complex plane $\mathbb C$ defined as

$$Tf(z)=\frac1{\pi}p.v.\int_{\mathbb C}\frac{f(\zeta)}{(z-\zeta)^{2}}d\zeta,$$
satisfies a certain sharp weighted estimate, and as a consequence they obtained borderline regularity properties for solutions of the Beltrami equation on $\mathbb C$ ($u_{\overline{z}} = \mu u_{z},$ where $\mu$ is a given function of $\|\mu\|_{L^{\infty}}<1$).

Let us now consider a more general class of operators than that of Calder\'on-Zygmund operators.  There are operators that for example, are not bounded in all Lebesgue spaces, $L^{p}$, $1<p<\infty$, and therefore we can not expect these operators to satisfy the assumptions of the previous theorem, since these assumptions imply that our operator is bounded in all $L^{p}$, $1<p<\infty$. Thus, it is natural to try and investigate what happens if we have an operator that is bounded in the weighted $L^{p}(w)$ spaces, but only for a subclass of the $A_{p}$ space. What kind of boundedness properties does this operator have? We will answer this question with the following theorem, which is the main result of this paper.

\begin{theorem}
\label{main}
Suppose that for some $1<p_{0}<\infty$, an operator $T$ in $\mathbb R^n$ satisfies the following inequality,

$$\|T\|_{L^{p_{0}}(w)\rightarrow L^{p_{0},\infty}(w)}\leq F([w]_{A_{p_{0}}}),$$
for all $A_{p_{0}}$ ``flat" weights, $w$, where $F$ is a positive increasing function. Then for all $p$ in a neighborhood around $p_{0}$, the operator $T$ is weakly bounded on $L^{p}(w)$ for all sufficiently ``flat" weights $w\in A_{p}$. In particular we have the estimate 

$$\|T\|_{L^{p}(w)\rightarrow L^{p,\infty}(w)}\leq J([w]_{A_{p}}),$$
for ``flat" $A_{p}$ weights w, where $J$ is positive function that depends on $F$, the dimension $n$, $p$, and $p_{0}$. 
\end{theorem}

A result of the same nature can be found in \cite{AM}, where the authors prove that boundedness on some $L^{p_{0}}(w)$, for nice weights $w$, imply boundedness for $L^{p}$ for $p$ close to $p_{0}$ (page 243). Their assumption requires the weight to belong to some Reverse H\"older classes and then they are able to extrapolate using the well-known Rubio De Francia iteration algorithm. Our result is very similar, and the main idea is the same, but the assumptions and the tools that are used for the proof of Theorem \ref{main} are not. 

\begin{remark}
\label{remark}
This theorem and Marcinkiewicz interpolation  imply that under some weak assumptions on the boundedness of the operator on $L^{p_{0}}$, we get that the operator is strongly bounded on $L^{p}$, where $p$ lies in some neighborhood around $p_{0}$. This is a rather unexpected result. Note that in order to use the Marcinkiewicz interpolation theorem we need to assume that our operator satisfies some kind of sub-linearity because in order to interpolate we have to be able to control the image of the sum of two functions under the action of $T$, by the sum of the images of the two functions under $T$.  
\end{remark}

\begin{remark}
\label{remark2}
In the proof of our main theorem we do not require any structure (sub-linearity for example) for the operator $T$ that appears in the statement of Theorem \ref{main}. The only property needed is that it is weakly bounded on the weighted $L^{p_{0}}$ space, as the same situation appears in the proof of the classical extrapolation theorem.  
\end{remark}
For the proof of this theorem we are going to need the following result proven in \cite{NV1}.

\begin{theorem}
\label{oldresult}
Suppose that $1<p<\infty$. There exists a constant $c_{n,p}>0$ that depends only on the dimension, $n,$ and on $p$ such that

$$\|M\|_{L^{p}(w)\rightarrow L^{p}(w)}\leq\|M\|_{L^{p}\rightarrow L^{p}}(1+c_{n,p}\sqrt{[w]_{A_{p}}-1}),$$
for all ``flat" weights $w\in A_{p}$, where $M$ is the non-centered Hardy-Littlewood maximal function defined as

$$Mf(x)=\sup_{x\in Q}\frac1{|Q|}\int_{Q}|f(x)|dx,$$
and $Q$ denotes a cube in $\mathbb R^n$.
\end{theorem}

This theorem holds for a much larger class of operators, see \cite{NV1} and a generalization in \cite{NV2}, where the authors defined a metric structure on the $A_{p}$ classes and studied continuity properties of weighted estimates with respect to this metric. 

In the next section, we present the proof of our main theorem, and in Section \ref{sec3}, we discuss some observations on the $L^{p}$ operator norm of the Maximal function, and remarks about the main theorem.

\end{section}

\begin{section}{the proof of the main result}
\label{sec2}
Our proof follows the same steps as the one given in \cite{DMP} (page 40), but there are some small changes that are needed in order to be able to extrapolate for a smaller class of weights. 

We begin by considering a fixed number $0<\epsilon<1$, and two positive functions $g\in L^{p}(w)$ and $h\in L^{p'}(w)$, where $w$ is an $A_{p}$ weight. We will fix the numbers $\epsilon$ and $1<p<\infty$ later. Let us define the operators

$$R_{\epsilon}g(x)=\sum_{k=0}^{\infty}\frac{M^{k}g(x)}{(1+\epsilon)^{k}\|M\|_{p,w}^{k}},$$
and

$$R^{'}_{\epsilon}h(x)=\sum_{k=0}^{\infty}\frac{(M^{'})^{k}h(x)}{(1+\epsilon)^{k}\|M^{'}\|_{p^{'}w}}.$$
where by $\|M\|_{p,w}$ we denote the $L^{p}(w)$ operator norm of $M$, and by $M^{k}$ we denote the composition of $M$ with itself, $k$ times, $M^{k}=M\circ\dots\circ M$. Also, $M^{0}$ denotes the identity operator. Similarly for the $M^{'}$ operator, that is the dual operator of $M$, defined as

$$M^{'}f(x)=\frac{M(fw)(x)}{w(x)}.$$ 
Since $w^{1-p'}\in A_{p'}$ we have that $M$ is bounded on $L^{p'}(w^{1-p'})$. Equivalently, $M^{'}$ is bounded on $L^{p'}(w)$. 

The Neumann series that was just defined converges in $L^{p}(w)$ and $L^{p'}(w)$, respectively and they define positive functions $R_{\epsilon}\in L^{p}(w)$ and $R^{'}_{\epsilon}\in L^{p'}(w)$. It is obvious that the function $R_{\epsilon}$ satisfies the inequality

\begin{equation}
\label{eq1}
g(x)\leq R_{\epsilon}g(x)
\end{equation}
for almost every $x\in\mathbb R^{n}$, since all the terms of the series are positive numbers. Moreover, we claim that

\begin{equation}
\label{eq2}
\|R_{\epsilon}g\|_{L^{p}(w)}\leq\frac{1+\epsilon}{\epsilon}\|g\|_{L^{p}(w)}.
\end{equation}
Indeed
\begin{eqnarray*}
\|R_{\epsilon}g\|_{L^{p}(w)}=\Big\|\sum_{k=0}^{+\infty}\frac{M^{k}g}{(1+\epsilon)^{k}\|M\|_{p,w}^{k}}\Big\|_{L^{p}(w)}&\leq&\sum_{k=0}^{+\infty}\frac{\|M^{k}g\|_{L^{p}(w)}}{(1+\epsilon)^{k}\|M\|^{k}_{p,w}}\\
&\leq&\sum_{k=0}^{+\infty}\frac{\|M^{k}\|_{p,w}\|g\|_{L^{p}(w)}}{(1+\epsilon)^{k}\|M\|^{k}_{p,w}}\\
&\leq&\frac{1+\epsilon}{\epsilon}\|g\|_{L^{p}(w)},
\end{eqnarray*}
which proves the claim.

The last important property of the $R_{\epsilon}g$ function is that it is an $A_{1}$ weight. A weight $w$ is said to belong to the $A_{1}$ class if there exists a positive constant, $c,$ such that 

$$Mw(x)\leq cw(x),$$
for almost every $x\in\mathbb R^n$. The smallest such $c$ is called the $A_{1}$ characteristic of $w$ and is denoted by $[w]_{A_{1}}$. The last important property of $R_{\epsilon}$ is the following,

\begin{equation}
\label{eq3}
[R_{\epsilon}g]_{A_{1}}\leq(1+\epsilon)\|M\|_{p,w}.
\end{equation}
This is easy to see since,

\begin{eqnarray*}
M(R_{\epsilon}g)=M\Big(\sum_{k=0}^{\infty}\frac{M^{k}g}{(1+\epsilon)^{k}\|M\|_{p,w}^{k}}\Big)&\leq&\sum_{k=0}^{\infty}\frac{M^{k+1}g}{(1+\epsilon)^{k}\|M\|_{p,w}^{k}}\\
&=&(1+\epsilon)\|M\|_{p,w}\sum_{k=0}^{\infty}\frac{M^{k+1}g}{(1+\epsilon)^{k+1}\|M\|_{p,w}^{k+1}}\\
&=&(1+\epsilon)\|M\|_{p,w}(R_{\epsilon}g-g)\\
&\leq&(1+\epsilon)\|M\|_{p,w}R_{\epsilon}g,
\end{eqnarray*}
since the function $g$ is chosen to be positive. Obviously, estimates like (\ref{eq1}) and (\ref{eq2}) will be satisfied by $R^{'}_{\epsilon}$ and $h$ in the place of $R_{\epsilon}$ and $g$. Namely,

\begin{equation}
\label{eqeq}
h(x)\leq R^{'}_{\epsilon}h(x)
\end{equation} 
and

\begin{equation}
\label{eqeqeq}
\|R^{'}_{\epsilon}h\|_{L^{p'}(w)}\leq\frac{1+\epsilon}{\epsilon}\|h\|_{L^{p'}(w)}.
\end{equation}
Also note that $M^{'}(R^{'}_{\epsilon}h)(x)\leq (1+\epsilon)\|M^{'}\|_{p',w}R^{'}_{\epsilon}h(x)$, and so 

\begin{equation}
\label{eq4}
[(R^{'}_{\epsilon}h)w]_{A_{1}}\leq(1+\epsilon)\|M^{'}\|_{p',w}.
\end{equation}

Now we will continue with the proof of Theorem \ref{main}. We will use strong type estimates for the proof, but surprisingly we do not lose anything by doing this. Indeed, suppose that we have an operator $G$ that is weakly bounded on $L^{p}(w)$. We fix $\lambda>0$ and for a function $f\in L^{p}(w)$, we define the set $E_{\lambda,f}=\{x\in\mathbb R^n:|Gf(x)|>\lambda\}$. Then we claim that the operator

$$T_{\lambda}f(x):=\lambda\chi_{E_{\lambda,f}}(x),$$
is strongly bounded on $L^{p}(w)$, and actually we have the uniform estimate 

$$\sup_{\lambda>0}\|T_{\lambda}\|_{L^{p}(w)\rightarrow L^{p}(w)}=\|G\|_{L^{p}(w)\rightarrow L^{p,\infty}(w)}.$$
To see this let us calculate the following norm

\begin{eqnarray*}
\|T_{\lambda}f\|_{L^{p}(w)}^{p}=\int_{\mathbb R^n}|T_{\lambda}f(x)|^{p}dx&=&\lambda^{p}w(\{x\in\mathbb R^{n}:|Gf(x)|>\lambda\})\\
&\leq&\|G\|_{L^{p}(w)\rightarrow L^{p,\infty}(w)}^{p}\|f\|_{L^{p}(w)}^{p}.
\end{eqnarray*}

\begin{remark}
\label{remark3}
In general, the operators $T_{\lambda}$ are not sub-linear even if we assume that the operator $G$ is linear. They satisfy the inequality

$$T_{\lambda}(f+g)\leq T_{\lambda}(2f)+T_{\lambda}(2g).$$
\end{remark}
Now our point is obvious and so let us proceed with the proof of the main theorem assuming that our operator is strongly bounded on $L^{p}(w)$.

We need to consider two different cases: 1) $p<p_{0}$ and 2) $p_{0}<p$. First consider the case $p<p_{0}$, but sufficiently close (we will make this more precise later).  Let $f\in L^{p}(w)$ and denote $g := \frac{|f|}{\|f\|_{L^{p}(w)}}$. Our goal is to estimate the quantity $\|Tf\|_{L^{p}(w)}$. We will use H\"older's inequality with respect to the $wdx$ measure, with $\alpha=(p_{0}-p)\frac{p}{p_{0}}$. Analytically,

\begin{eqnarray*}
\|Tf\|_{L^{p}(w)}&=&\Big(\int_{\mathbb R^{n}}|Tf(x)|^{p}(R_{\epsilon}g(x))^{-\alpha}(R_{\epsilon}g(x))^{\alpha}w(x)dx\Big)^{\frac1{p}}\\
&\leq&\Big(\int_{\mathbb R^n}|Tf(x)|^{p_{0}}(R_{\epsilon}g(x))^{\frac{-\alpha p_{0}}{p}}w(x)dx\Big)^{\frac1{p_{0}}}\Big(\int_{\mathbb R^n}(R_{\epsilon}g(x))^{\alpha(\frac{p_{0}}{p})'}w(x)dx\Big)^{\frac1{p(\frac{p_{0}}{p})'}}\\
&=&\Big(\int_{\mathbb R^n}|Tf(x)|^{p_{0}}(R_{\epsilon}g(x))^{p-p_{0}}w(x)dx\Big)^{\frac1{p_{0}}}\Big(\int_{\mathbb R^n}(R_{\epsilon}g(x))^{p}w(x)dx\Big)^{\frac{p_{0}-p}{pp_{0}}}.
\end{eqnarray*}
Denote $W(x):=(R_{\epsilon}g(x))^{p-p_{0}}w(x)$. We claim that $W\in A_{p_{0}}$ and furthermore that $W$ is ``flat", for $p$ close enough to $p_{0}$. Indeed,

$$\frac1{|Q|}\int_{Q}(R_{\epsilon}g(x))^{-(p_{0}-p)}w(x)dx\leq [R_{\epsilon}g]_{A_{1}}^{p_{0}-p}\Big(\frac1{|Q|}\int_{Q}w(x)dx\Big)\Big(\frac1{|Q|}\int_{Q}R_{\epsilon}g(x)dx\Big)^{-(p_{0}-p)},$$
and

\begin{eqnarray*}
\frac1{|Q|}\int_{Q}W(x)^{1-p_{0}'}dx&=& \frac1{|Q|}\int_{Q}(R_{\epsilon}g(x))^{\frac{p_{0}-p}{p_{0}-1}}w(x)^{1-p_{0}'}dx\\
&\leq&\Big(\frac1{|Q|}\int_{Q}R_{\epsilon}g(x)dx\Big)^{\frac{p_{0}-p}{p_{0}-1}}\Big(\frac1{|Q|}\int_{Q}w(x)^{1-p'}dx\Big)^{\frac{p-1}{p_{0}-1}}.
\end{eqnarray*}
Therefore,

\begin{eqnarray*}
\Big(\frac1{|Q|}\int_{Q}W(x)dx\Big)\Big(\frac1{|Q|}\int_{Q}W(x)^{1-p_{0}'}dx\Big)^{p_{0}-1}&\leq&[R_{\epsilon}g]_{A_{1}}^{p_{0}-p}\Big(\frac1{|Q|}\int_{Q}w(x)dx\Big)\\
&\cdot&\Big(\frac1{|Q|}\int_{Q}w(x)^{1-p'}dx\Big)^{p-1},
\end{eqnarray*}
which means that (we used inequality (\ref{eq3}))

$$[W]_{A_{p_{0}}}\leq(1+\epsilon)^{p_{0}-p}\|M\|_{p,w}^{p_{0}-p}[w]_{A_{p}}.$$
For $p$ close to $p_{0}$ and for $w$ ``flat" in $A_{p}$ we have that $W$ is ``flat" in $A_{p_{0}}$ (since we can control $\|M\|_{p,w}$ by Theorem \ref{oldresult}). Now we estimate $\|Tf\|_{L^{p}(w)}$ as

\begin{eqnarray*}
\|Tf\|_{L^{p}(w)}&\leq&\Big(\int_{\mathbb R^n}|Tf(x)|^{p_{0}}W(x)dx\Big)^{\frac1{p_{0}}}\Big(\frac{1+\epsilon}{\epsilon}\Big)^{\frac{p_{0}-p}{p_{0}}}\|g\|_{L^{p}(w)}^{\frac{p_{0}-p}{p_{0}}}\\
&\leq&\Big(\int_{\mathbb R^n}|f(x)|^{p_{0}}(R_{\epsilon}g(x))^{p-p_{0}}w(x)dx\Big)^{\frac1{p_{0}}}\Big(\frac{1+\epsilon}{\epsilon}\Big)^{\frac{p_{0}-p}{p_{0}}}F([W]_{A_{p_{0}}})\\
&\leq&\Big(\int_{\mathbb R^n}|f(x)|^{p_{0}}g(x)^{p-p_{0}}w(x)dx\Big)^{\frac1{p_{0}}}\Big(\frac{1+\epsilon}{\epsilon}\Big)^{\frac{p_{0}-p}{p_{0}}}F([W]_{A_{p_{0}}})\\
&=&\Big(\int_{\mathbb R^n}|f(x)|^{p_{0}}\frac{|f(x)|^{p-p_{0}}}{\|f\|_{L^{p}(w)}^{p-p_{0}}}w(x)dx\Big)^{\frac1{p_{0}}}\Big(\frac{1+\epsilon}{\epsilon}\Big)^{\frac{p_{0}-p}{p_{0}}}F([W]_{A_{p_{0}}})\\
&=&\Big(\int_{\mathbb R^n}\frac{|f(x)|^{p}}{\|f\|_{L^{p}(w)}^{p-p_{0}}}w(x)dx\Big)^{\frac1{p_{0}}}\Big(\frac{1+\epsilon}{\epsilon}\Big)^{\frac{p_{0}-p}{p_{0}}}F([W]_{A_{p_{0}}})\\
&=&\frac{\|f\|_{L^{p}(w)}^{\frac{p}{p_{0}}}}{\|f\|_{L^{p}(w)}^{\frac{p-p_{0}}{p_{0}}}}\Big(\frac{1+\epsilon}{\epsilon}\Big)^{\frac{p_{0}-p}{p_{0}}}F([W]_{A_{p_{0}}})\\
&\leq&\Big(\frac{1+\epsilon}{\epsilon}\Big)^{\frac{p_{0}-p}{p_{0}}} F((1+\epsilon)^{p_{0}-p}\|M\|_{p,w}^{p_{0}-p}[w]_{A_{p}})\|f\|_{L^{p}(w)}.
\end{eqnarray*} 
which means that we are done. Note that here we used inequality (\ref{eq1}) and the fact that the exponent $p-p_{0}$ is negative to substitute the function $R_{\epsilon}g$ by $g$ and then use that $g=\frac{|f|}{\|f\|_{L^{p}(w)}}$.

Finally, we need to deal with the case for which $p_{0}<p$. Again we begin with a function $f\in L^{p}(w)$. By duality there exists a non-negative function $h\in L^{p'}(w)$, $\|h\|_{L^{p'}(w)}=1$, with the property

$$\|Tf\|_{L^{p}(w)}=\int_{\mathbb R^n}Tf(x)h(x)w(x)dx.$$
Our claim is that the weight $W(x)=(R_{\epsilon}h(x))^{\frac{p-p_{0}}{p-1}}w(x)\in A_{p_{0}}$, with characteristic very close to $1$, if $p$ is sufficiently close to $p_{0}$, with $w$ ``flat" in $A_{p}$. Indeed, we set $q=\frac{p-1}{p-p_{0}}>1$, $q'=\frac{p-1}{p_{0}-1}$, and by H\"older's inequality we obtain 

\begin{eqnarray*}
\frac{1}{|Q|}\int_{Q}W(x)dx&=&\frac1{|Q|}\int_{Q}(R_{\epsilon}h(x))^{\frac{p-p_{0}}{p-1}}w(x)dx\\
&\leq&\Big(\frac1{|Q|}\int_{Q}(R_{\epsilon}h(x))w(x)dx\Big)^{\frac{p-p_{0}}{p-1}}\Big(\frac1{|Q|}\int_{Q}w(x)dx\Big)^{\frac{p_{0}-1}{p-1}},
\end{eqnarray*}  
and 

\begin{eqnarray*}
\Big(\frac1{|Q|}\int_{Q}W(x)^{1-p_{0}'}dx\Big)^{p_{0}-1}&=&\Big(\frac1{|Q|}\int_{Q}(R_{\epsilon}h(x))^{-\frac{p-p_{0}}{(p-1)(p_{0}-1)}}w(x)^{1-p_{0}'}dx\Big)^{p_{0}-1}\\
&\leq&[(R^{'}_{\epsilon}h)w]_{A_{1}}^{\frac{p-p_{0}}{p-1}}\Big(\frac1{|Q|}\int_{Q}(R^{'}_{\epsilon}h(x))w(x)dx\Big)^{-\frac{p-p_{0}}{p-1}}\\
&\cdot&\Big(\frac1{|Q|}\int_{Q}w(x)^{\frac{p-p_{0}}{(p-1)(p_{0}-1)}}w(x)^{1-p_{0}'}dx\Big)^{p_{0}-1}\\
&=&[(R^{'}_{\epsilon}h)w]_{A_{1}}^{\frac{p-p_{0}}{p-1}}\Big(\frac1{|Q|}\int_{Q}(R^{'}_{\epsilon}h(x))w(x)dx\Big)^{-\frac{p-p_{0}}{p-1}}\\
&\cdot&\Big(\frac1{|Q|}\int_{Q}w(x)^{1-p'}dx\Big)^{p_{0}-1}.
\end{eqnarray*}
Gathering everything together we have (we used inequality (\ref{eq4}))

$$[W]_{A_{p_{0}}}\leq(1+\epsilon)^{\frac{p-p_{0}}{p-1}}\|M\|_{L^{p'}(w^{1-p'})}^{\frac{p-p_{0}}{p-1}}[w]_{A_{p}}^{\frac{p_{0}-1}{p-1}}.$$
We see that for $p\approx p_{0}$, and for $w$ ``flat" in $A_{p}$ we have that $W$ is ``flat" in $A_{p_{0}}$ (we control the norm of the Maximal function by Theorem \ref{oldresult}), which means that we can use our boundedness assumption for the operator $T$ and of course H\"older's inequality. Therefore,

\begin{eqnarray*}
\int_{\mathbb R^n} Tf(x)h(x)w(x)dx&\leq&\int_{\mathbb R^n}Tf(x) R^{'}_{\epsilon}h(x)^{\frac{p-p_{0}}{p_{0}(p-1)}}h(x)^{\frac{(p_{0}-1)p}{p_{0}(p-1)}}w(x)dx\\
&\leq&\Big(\int_{\mathbb R^n}Tf(x)^{p_{0}}W(x)dx\Big)^{\frac1{p_{0}}}\Big(\int_{\mathbb R^n}h(x)^{p'}w(x)dx\Big)^{\frac1{p_{0}^{'}}}\\
&=&\Big(\int_{\mathbb R^n}Tf(x)^{p_{0}}W(x)dx\Big)^{\frac1{p_{0}}}\\
&\leq&F([W]_{A_{p_{0}}})\Big(\int_{\mathbb R^n}f(x)^{p_{0}}R^{'}_{\epsilon}h(x)^{\frac{p-p_{0}}{p-1}}w(x)dx\Big)^{\frac1{p_{0}}}\\
&\leq&F\Big((1+\epsilon)^{\frac{p-p_{0}}{p-1}}\|M\|_{L^{p'}(w^{1-p'})}^{\frac{p-p_{0}}{p-1}}[w]_{A_{p}}^{\frac{p_{0}-1}{p-1}}\Big)\|f\|_{L^{p}(w)}\|R^{'}_{\epsilon}h\|_{L^{p'}(w)}^{\frac{p-p_{0}}{p_{0}(p-1)}}\\
&\leq&\Big(\frac{1+\epsilon}{\epsilon}\Big)^{\frac{p-p_{0}}{p_{0}(p-1)}}F\Big((1+\epsilon)^{\frac{p-p_{0}}{p-1}}\|M\|_{L^{p'}(w^{1-p'})}^{\frac{p-p_{0}}{p-1}}[w]_{A_{p}}^{\frac{p_{0}-1}{p-1}}\Big)\|f\|_{L^{p}(w)},
\end{eqnarray*} 
which means we are done (we used inequality (\ref{eqeqeq})).

\end{section}

\begin{section}{remarks on the main theorem and the maximal function}
\label{sec3}
In this section we will discuss the asymptotic behavior of the operator norm of the Maximal operator in $L^{p}$, and then make some observations about the main Theorem \ref{main}. 

Let us mention that it is natural to be interested in how the operator norm of the Maximal function behaves with respect to $p$, since the Maximal operator controls so many other operators of Harmonic Analysis. In \cite{GM} the operator norm of $M$ was found for all $L^{p}$ spaces, $1<p<+\infty$, in dimension $1$. The operator norm of the maximal operator is not known in any of the $L^{p}$ spaces for dimension $n\geq2$, with $1<p<\infty$. It is well known that 

$$\lim_{p\to+\infty}\|M\|_{L^{p}\rightarrow L^{p}}=1,$$
and we refer the reader to \cite{GM} for more results and the reference therein.  

Once again, notice that Theorem \ref{main} under some weak assumptions on the weak boundedness of the operator on $L^{p_{0}}$, implies boundedness on $L^{q}$ for all $q$ that are sufficiently close to $p_{0}$. This is a rather unexpected result (see Remark \ref{remark}).  Moreover, the proof of this theorem follows the same lines as the proof of the classical extrapolation, Theorem \ref{result1}, but in order to get our result, we had to make some small changes to ensure that our weight $W$ that appears naturally in the calculations (see section \ref{sec2}), is sufficiently ``flat". 

Let us note that from the proof of Theorem \ref{main} follows that the diameter of the neighborhood around $p_{0}$ in which the operator is weakly bounded (and under the sub-linearity assumption, strongly bounded (see Remark \ref{remark})) depends only on the quantitative behavior of the operator norm of the Maximal function and the initial assumption on the ``flatness" of the $L^{p_{0}}(w)$ estimate. To see this, choose $p$ so close to $p_{0}$ such that $\|M\|$ raised to $\max\{p-p_{0},p_{0}-p\}$ is close to $1$. Then, choose the weight $w$ to be sufficiently ``flat" in $A_{p}$, make the parameter $\epsilon>0$ that appears in the proof sufficiently small, and  and we are done. We have pointed out that $\lim_{p\to+\infty}\|M\|_{L^{p}\rightarrow L^{p}}=1$, which means that we expect the neighborhood to become larger on the left of $p_{0}$, as we consider larger values of the number $p_{0}$, and larger on the right of $p_{0}$ as we consider values of the number $p_{0}$ closer to $1$. This follows from the fact that in the case where $p<p_{0}$ the norm $\|M\|_{L^{p}\rightarrow L^{p}}$ appears in the estimates, and in the case $p_{0}<p$ the norm $\|M\|_{L^{p'}\rightarrow L^{p'}}$ comes out naturally from the estimates.

\end{section}

Nicholas Boros, borosnic@msu.edu\newline

Nikolaos Pattakos, pattakos@msu.edu\newline
 
Alexander Volberg, volberg@math.msu.edu

\end{document}